\documentclass[10pt]{amsart}
\usepackage{amsmath}
\usepackage{amsfonts}
\usepackage{latexsym}
\usepackage{amssymb}
\usepackage{enumerate}
\usepackage[dvips]{graphics}

\newcommand{\Z}{{\mathbf Z}}

\newcommand{\R}{{\mathbf R}}
\newcommand{\C}{{\mathbf C}}

\newcommand{\rk}{{\rm rk}}
\newcommand{\comment}[1]{}

\def\V{\mathcal V}

\def\ker{{\rm Ker }}

\def\R{\mathbf R}

\newtheorem{theorem}{Theorem}
\newtheorem{proposition}{Proposition}

\newtheorem{corollary}[proposition]{Corollary}

\begin{document}
\title{Novikov-Betti numbers and the fundamental group}


\author[M. ~Farber and D. ~Sch\"utz]{M.~Farber and D. ~Sch\"utz}
\address{Department of Mathematics, University of Durham, Durham DH1 3LE, UK}
\email{Michael.Farber@durham.ac.uk} \email{Dirk.Schuetz@durham.ac.uk}

\maketitle

\begin{theorem}\label{thm1}
Let $X$ be a connected finite polyhedron and let $\xi\in H^1(X;\R)$ be a
nonzero cohomology class. If the first Novikov-Betti number $b_1(\xi)$ is
nonzero, $b_1(\xi)>0$, then $\pi_1(X)$ contains a nonabelian free subgroup.
\end{theorem}

This result may appear striking as the Novikov-Betti numbers carry \lq\lq
abe\-lian\rq\rq \, information about $X$. We refer to \cite{noviko},
\cite{farber} for the definition of the Novikov-Betti numbers; an explicit
definition will also be given below in the proof of Theorem \ref{thm1}.

An alternative description of $b_i(\xi)$ uses homology of complex flat line
bundles. Consider the variety $\V_\xi$ of all complex flat line bundles $L$
over $X$ having the following property: $L$ has trivial monodromy along any
loop $\gamma$ in $X$ assuming that $\langle \xi, [\gamma]\rangle =0$. It is
easy to see that (a) $\V_\xi$ is an algebraic variety isomorphic to
$(\C^\ast)^r$ for some integer $r$ and (b) the dimension $\dim H_i(X;L)$ is
independent of $L$ assuming that $L\in \V_\xi$ is generic, see \cite{farber},
Theorem 1.50. The number $\dim H_i(X;L)$ for a generic $L\in \V_\xi$ coincides
with the Novikov-Betti number $b_i(\xi)$.

The proof of Theorem \ref{thm1} given below is based on the results of R.
Bieri, W. Neumann, R. Strebel \cite{BNS} and J.-Cl. Sikorav \cite{Si1}.

Theorem \ref{thm1} implies the following vanishing result:

\begin{corollary}\label{cor1}
Let $X$ be a connected finite polyhedron having an amenable fundamental group.
Then the first Novikov-Betti number vanishes $b_1(\xi)=0$ for any $\xi\not=0\in
H^1(X;\R)$.
\end{corollary}

\noindent Corollary \ref{cor1} follows from Theorem \ref{thm1} as an amenable
group contains no nonabelian free subgroups. As another simple corollary of
Theorem \ref{thm1} we mention the next statement:

\begin{corollary}\label{cor2}
Assume that $X$ is a connected finite two-dimensional polyhedron. If the Euler
characteristic of $X$ is negative $\chi(X)<0$ then $\pi_1(X)$ contains a
nonabelian free subgroup.
\end{corollary}

\begin{proof} The first Betti number $b_1(X)$ is positive as follows from
$\chi(X)<0$. Hence there exists a nonzero cohomology class $\xi\in H^1(X;\R)$.
Then \begin{eqnarray}\label{neg} \chi(X)= b_2(\xi)-b_1(\xi),\end{eqnarray} see
Proposition 1.40 in \cite{farber}. Here we use that $b_0(\xi)=0$ for
$\xi\not=0$, see Corollary 1.33 in \cite{farber}. The inequality $\chi(X)<0$
together with (\ref{neg}) imply that $b_1(\xi)>0$. Theorem \ref{thm1} now
states that $\pi_1(X)$ contains a nonabelian free subgroup.
\end{proof}
\noindent Corollary \ref{cor2} is known, it was obtain by N.S. Romanovskii, see
\cite{Ro}. It can be equivalently expressed algebraically as follows:

\begin{corollary}
Any discrete group $G$ of deficiency greater than 1 contains a nonabelian free
subgroup.
\end{corollary}

The authors thank Jonathan Hillman for referring to \cite{Ro}.


\begin{proof}[Proof of Theorem \ref{thm1}]
Consider the following diagram of rings and ring homomorphisms:
\begin{eqnarray}\label{diagram}
\begin{array}{lcr}
& \Lambda & \\
\,\, \alpha \swarrow && \searrow \rho\, \, \, \\
S_\xi & \stackrel \beta \longrightarrow & N_\xi
\end{array}
\end{eqnarray}

Here $\Lambda=\Z[\pi]$ is the group ring of the fundamental group of $X$ which
we denote $\pi=\pi_1(X,x_0)$.

The ring $N_\xi$ is the {\it Novikov ring} which is defined as follows. View
the class $\xi$ as a group homomorphism $\xi: \pi\to \Z$ and let $H$ denotes
the factorgroup $\pi/\ker(\xi)$. It is a finitely generated free abelian group.
The usual group ring $\Z[H]$ consists of finite sums of the form $\sum a_j h_j$
with $a_j\in \Z$ and $h_j\in H$; it coincides with the Laurent polynomial ring
in $\rk H$ variables. The Novikov ring $N_\xi$ is a completion of $\Z[H]$; its
elements are infinite sums $\sum_{j=1}^\infty a_j h_j$ with $a_j\in \Z$ and
$h_j\in H$ such that the sequence of evaluation $\xi(h_j)$ tends to $-\infty$.
In other words, $N_\xi$ is the ring of Laurent power series in $\rk H$
variables where the terms of the series go to infinity in the direction
specified by the class $\xi$.

The ring $S_\xi$ is the {\it Novikov - Sikorav} completion of the group ring
$\Lambda=\Z[\pi]$; it was originally introduced by J.-Cl. Sikorav \cite{Si1}.
Elements of $S_\xi$ are infinite sums of the form $\sum_{j=1}^\infty a_j g_j$
where $a_j\in \Z$, $g_j\in \pi$ and $\xi(g_j)\to -\infty$.

Diagram (\ref{diagram}) includes the obvious rings homomorphisms $\alpha$,
$\beta$, $\rho$ and is commutative.

Let $X$ be a finite connected polyhedron with fundamental group $\pi$. Consider
the universal covering $\tilde X\to X$. The cellular chain complex
$C_\ast(\tilde X)$ is a complex of finitely generated free left
$\Lambda$-modules. Since $N_\xi$ is a commutative ring and a principal ideal
domain (see \cite{farber}, Lemma 1.10) the homology $H_i(N_\xi\otimes_\Lambda
C_\ast(\tilde X))$ is a finitely generated $N_\xi$-module and its rank (over
$N_\xi$) equals the Novikov-Betti number $b_i(\xi)$, see \cite{farber}, \S 1.5.

Our first goal is to show that the assumption $b_1(\xi)\not=0$ implies
\begin{eqnarray}\label{first} H_1(S_\xi\otimes_\Lambda C_\ast(\tilde X))\not=0.\end{eqnarray}
Indeed, note that $H_0(S_\xi\otimes_\Lambda C_\ast(\tilde X))=0$ (as we assume
that $\xi\not=0$). Hence, using the K\"unneth spectral sequence we find that
\begin{eqnarray}\label{tensor}
H_1(N_\xi\otimes_\Lambda C_\ast(\tilde X))
 =
 H_1(N_\xi\otimes_{S_\xi}(S_\xi\otimes_\Lambda
C_\ast(\tilde X))) \nonumber \\= N_\xi\otimes_{S_\xi} H_1(S_\xi\otimes_\Lambda
C_\ast(\tilde X)).
\end{eqnarray}
Hence $b_1(\xi)\not=0$ implies the nonvanishing of $H_1(N_\xi\otimes_\Lambda
C_\ast(\tilde X))$ which via (\ref{tensor}) gives (\ref{first}).

Note that $b_1(\xi)=b_1(-\xi)$ (see Corollary 1.31 in \cite{farber}) and hence
the nonvanishing of the Novikov-Betti number $b_1(\xi)$ implies the
nonvanishing of the Novikov - Sikorav homology $H_1(S_\xi\otimes_\Lambda
C_\ast(\tilde X))$ for both cohomology classes $\xi$ and $-\xi$.

R. Bieri, W. Neumann and R. Strebel \cite{BNS} defined a geometric invariant of
discrete groups. It can be viewed as a subset $\Sigma$ of the space of nonzero
cohomology classes $\Sigma\subset H^1(X;\R)-\{0\}$. A theorem of Jean-Claude
Sikorav (see \cite{Si1}, page 86 or \cite{bieri}) states that for a nonzero
cohomology class $\xi\in H^1(X;\R)$ the following conditions are equivalent:
(a) $H_1(S_\xi\otimes_{\Lambda}C_\ast(\tilde X))=0$ and (b) $-\xi\in \Sigma$.
Hence, as explained above, $b_1(\xi)\not=0$ implies that $\xi\not\in \Sigma$
and $-\xi\not\in \Sigma$. Now we apply Theorem C from \cite{BNS} which states
that the union of $\Sigma$ and $-\Sigma$ equals $H^1(X;\R)-\{0\}$ assuming that
$\pi_1(X)$ has no non-abelian free subgroups. Since in our case neither $\xi\in
\Sigma$ nor $-\xi \in \Sigma$ we conclude that $\pi_1(X)$ must contain a
non-abelian free subgroup. 
\end{proof}

Next we mention an example showing that a space with nonvanishing Novikov
torsion number $q_1(\xi)$ may have an amenable fundamental group. In other
words, Theorem \ref{thm1} becomes false if one replaces the assumption
$b_1(\xi)\not=0$ by the assumption $q_1(\xi)\not=0$. Example 1.49 in
\cite{farber} gives a two-dimensional polyhedron with $b_1(\xi)=0$ and
$q_1(\xi)=1$ for some nonzero $\xi\in H^1(X;\R)$. The fundamental group of $X$
is a Baumslag-Solitar group $G=|a, b; aba^{-1}=b^2|$. The later group appears
in the exact sequence $0\to \Z[\frac{1}{2}]\to G\to \Z\to 0$ (see
\cite{farber}, page 30); hence $G$ is amenable.

\bibliographystyle{amsalpha}

\end{document}